\documentclass[11pt]{article}
\usepackage{epic,latexsym,amssymb}
\usepackage{color}
\usepackage{amsfonts,amsmath,amssymb}
\usepackage{amsmath,amsthm}
\usepackage{comment}
\usepackage{pstricks,pst-node,pst-plot} 
\usepackage{tikz}

\newtheorem{thm}{Theorem}[section]
\newtheorem{lem}[thm]{Lemma}
\newtheorem{observation}[thm]{Observation}
\newtheorem{cor}[thm]{Corollary}

\newtheorem{ex}[thm]{Example}
\newtheorem{claim}{Claim}
\newtheorem{defn}{Definition}

\newcommand{\sta}{{\rm sta}}
\newcommand{\sdt}{{\rm sd_{\gamma_t}}}
\newcommand{\cT}{\mathcal{T}}
\newcommand{\cO}{\mathcal{O}}
\newcommand{\cP}{\mathcal{P}}

\newcommand{\2}{ \vspace{0.2cm} }
\newcommand{\1}{ \vspace{0.1cm} }
\newcommand{\smallqed}{{\tiny ($\Box$)}}
\newcommand{\modo}{{\rm mod \,}}

\let\oldenumerate\enumerate
\renewcommand{\enumerate}{
  \oldenumerate
  \setlength{\itemsep}{0pt}
  \setlength{\parskip}{0pt}
  \setlength{\parsep}{0pt}
}

\begin{document}

\title{On total domination subdivision numbers of trees}

\author{$^1$Michael A. Henning\thanks{Research supported in part by the South African National Research Foundation under grant number 129265 and the University of Johannesburg
} \,
and $^2$Jerzy Topp \\
\\
$^1$Department of Mathematics and Applied Mathematics\\
University of Johannesburg \\
Auckland Park, 2006, South Africa \\
\small \tt Email: mahenning@uj.ac.za \\
\\
$^2$Institute of Applied Informatics\\
University of Applied Sciences\\
82-300 Elbl\c{a}g, Poland \\
\small \tt Email: j.topp@ans-elblag.pl}

\date{}
\maketitle

\begin{abstract} A set $S$ of vertices in a graph $G$ is a total
dominating set of $G$ if every vertex is adjacent to a vertex in $S$.
The total domination number $\gamma_t(G)$ is the minimum cardinality of
a total dominating set of $G$. The total domination subdivision
number $\mbox{sd}_{\gamma_t}(G)$ of a graph $G$ is the minimum number
of edges that must be subdivided (where each edge in $G$ can be
subdivided at most once) in order to increase the total domination
number. Haynes et al. (Discrete Math. 286 (2004) 195--202) have given
a constructive characterization of trees whose total domination
subdivision number is~$3$. In this paper, we give new characterizations
of trees whose total domination subdivision number is 3.
\end{abstract}

{\small \textbf{Keywords:} Trees; Total domination number; Total domination subdivision number} \\
\indent {\small \textbf{AMS subject classification:} 05C69}

\section{Introduction}

For graph theory notation and terminology, we generally follow~\cite{HaHeHe-23}.  Specifically, let $G$ be a graph with vertex set $V(G)$ and edge set $E(G)$, and of order $n(G) = |V(G)|$ and size $m(G) = |E(G)|$. Two vertices $u$ and $v$ of $G$ are \emph{adjacent} if $uv\in E(G)$, and are called \emph{neighbors}. The \emph{open neighborhood} $N_G(v)$ of a vertex $v$ in $G$ is the set of neighbors of $v$, while the \emph{closed neighborhood} of $v$ is the set $N_G[v] = \{v\} \cup N_G(v)$. In general, for a subset $X \subseteq V(G)$, its \emph{open neighborhood} is the set $N_G(X) = \bigcup_{v \in X} N_G(v)$, and its \emph{closed neighborhood} is the set $N_G[X] = N_G(X) \cup X$.

The \emph{degree} of a vertex $v$ in $G$ is the number of neighbors of $v$ in $G$, and is denoted by $\deg_G(v)$, and so $\deg_G(v) = |N_G(v)|$. An \emph{isolated vertex} in $G$ is a vertex of degree zero. A~graph without any isolated vertex is called an \emph{isolate-free graph}. A vertex of degree~$1$ is called a \emph{leaf}, and its (unique) neighbor is called a \emph{support vertex}. The edge incident with a support vertex  and a leaf neighbor of the support vertex is called a \emph{pendant edge}. A \emph{strong support vertex} is a vertex with at least two leaf neighbors, and a \emph{weak support vertex} is a vertex with exactly one leaf neighbor. The set of leaves and the set of support vertices of $G$ are denoted by $L(G)$ and $S(G)$, respectively.

A graph $G$ is \emph{connected} if there is a $(u,v)$-path in $G$ joining every two vertices $u$ and $v$ in $G$. The \emph{distance} between two vertices $u$ and $v$ in a connected graph $G$, denoted by $d_G(u,v)$, is the minimum length among all $(u,v)$-paths in $G$. If $X$ and $Y$ are subsets of vertices of $G$, then the \emph{distance} $d_G(X,Y)$ between $X$ and $Y$ in $G$ is the minimum distance $d_G(x,y)$ among all pairs of vertices where $x \in X$ and $y \in Y$. The \emph{distance} $d_G(e_1,e_2)$ between two edges $e_1 = u_1v_1$ and $e_2 = u_2v_2$ of $G$ is the distance between the sets $\{u_1,v_1\}$ and $\{u_2,v_2\}$. The  \emph{distance} $d_G(e,F)$ between an edge $e$ and a subset $F$ of edges in $G$ is the minimum distance $d_G(e,f)$ between the edge $e$ and all edges $f \in F$. If $k$ is a positive integer and $u$ is a vertex in $G$, then the $k$-\emph{neighborhood} of $u$, denoted by $N_G^k(u)$, is the set of vertices at distance $k$ from $u$, that is, $N_G^k(u)=\{x\in V(G)\colon d_G(u,x)=k\}$.

We use $P_n$, $C_n$, and $K_n$ to denote a \emph{path}, a \emph{cycle}, and a \emph{complete graph}, respectively, on $n$ vertices. The \emph{complete bipartite graph} $K_{r,s}$ is a bipartite graph with partite sets $X$ and $Y$, where $|X| = r$, $|Y| = s$, and every vertex in $X$ is adjacent to every vertex in $Y$. A \emph{star} is a tree with at most one vertex that is not a leaf; that is, stars consist of complete bipartite graphs $K_{1,s}$ for $s \ge 1$ along with the trivial graph $K_1$. For $k \ge 1$ an integer, we let $[k]$ denote the set $\{1,\ldots,k\}$ and we let $[k]_0 = [k] \cup \{0\} = \{0,1,\ldots,k\}$.

A \emph{total dominating set}, abbreviated TD-set, of a graph $G$ is a set $S$ of vertices of $G$ such that every vertex has a neighbor in $S$, and so $N_G(S) = V(G)$. The \emph{total domination number} of $G$, denoted by $\gamma_t(G)$, is the minimum cardinality of a TD-set in $G$. For fundamentals on total domination theory in graphs we refer the reader to the authors' book~\cite{HeYe-13}, and to the so-called ``domination books''~\cite{HaHeHe-20,HaHeHe-21,HaHeHe-23}.

The \emph{total domination subdivision number} $\sdt(G)$ of a graph $G$ is the minimum number of edges that must be subdivided (where each edge in $G$ can be subdivided at most once) in order to increase the total domination number. The total domination subdivision number was introduced by Haynes et al.~\cite{Haynes-Hedetniemi-van der Merwe} and is now well studied in the literature, see~\cite{Favaron-Karami-Khoeilar-Sheikholeslami-1}--\cite{Favaron-Karami-Sheikholeslami-2},~\cite{Haynes-Hedetniemi-van der Merwe}--\cite{Haynes-Henning-Hopkins},~\cite{Jeyanthi-Hemalatha-Davvaz}--\cite{Soltankhah}, to mention but a few papers on the topic. Haynes et al.~\cite{Haynes-Hedetniemi-van der Merwe} have also shown that the total domination subdivision number of a~graph can be arbitrarily large, but the total domination subdivision number of a tree is either~$1$,~$2$, or~$3$, and so trees can be classified as Class~$1$, Class~$2$, or Class~$3$ depending on whether their total domination subdivision number is~$1$,~$2$, or~$3$, respectively.

A constructive characterization of all trees in Class~$3$ (Class~$1$, respectively) has been provided by Haynes, Henning, and Hopkins in~\cite{Haynes-Henning-Hopkins} (Karami, Khodkar, Khoeilar, and Sheikholeslami in~\cite{Karami-Khodkar-Khoeilar-Sheikholeslami}, respectively). In this paper, inspired by results of Haynes, Henning, and Hopkins~\cite{Haynes-Henning-Hopkins}, we continue their studies of trees in Class~$3$ and, in particular, we provide a different new characterization of the trees that belong to this class.

\subsection{Known results}

In this section, we present the constructive characterization of trees in Class~$3$ presented in~\cite{Haynes-Henning-Hopkins}. For this purpose the authors of~\cite{Haynes-Henning-Hopkins} describe a procedure to build a family $\cT$ of labeled trees that are of Class~$3$ as follows. We assign to each vertex $v$ a \emph{label}, also called its \emph{status}, and denoted by $\sta(v)$.

\begin{defn}
\label{defn:cT}
Let $\cT$ be the family of labeled trees that:
\begin{itemize}
\item[$(1)$] contains a path $P_6$ where the two leaves have status C, the two support vertices have status B, and the two central vertices have status A; and
\item[$(2)$] is closed under the two operations $\cT_1$ and $\cT_2$, which extend the labeled tree $T$ belonging to $\cT$ by attaching a tree to the vertex $y \in V(T)$.
\end{itemize}
\begin{itemize}
\item {\bf Operation} $\cT_1$. Assume $\sta(y) = A$. Then add a path $xwv$ and the edge $xy$. Let $\sta(x) = A$, $\sta(w) = B$, and $\sta(v) = C$.
\item  {\bf Operation} $\cT_2$. Assume $\sta(y)\in \{B,C\}$. Then add a path $xwvu$ and the edge $xy$. Let $\sta(x) = \sta(w) = A$, $\sta(v) = B$, and $\sta(u) = C$.
\end{itemize}
\end{defn}

The operations $\cT_1$ and $\cT_2$, and  a labeled tree belonging to the family $\cT$ are illustrated in Fig. \ref{tree-1}. The \emph{underlying tree} of a labeled tree $T$ is the tree obtained from $T$ by removing the vertex labels.

\begin{figure}[htb] \begin{center} 
\begin{tikzpicture}[scale=0.8,style=thick,x=0.1cm,y=0.1cm]
\def\vr{2.5pt} 
\draw (-5,0) node {\begin{tikzpicture}[scale=1,style=thick,x=0.1cm,y=0.1cm] \def\vr{2.0pt}
\draw (20,10)--(60,10); \draw (20,35)--(50,35);
\draw (20,10) [fill=black] circle (\vr);
\draw (30,10) [fill=black] circle (\vr);
\draw (40,10) [fill=black] circle (\vr);
\draw (50,10) [fill=black] circle (\vr);
\draw (60,10) [fill=black] circle (\vr);
\draw (20,35) [fill=black] circle (\vr);
\draw (30,35) [fill=black] circle (\vr);
\draw (40,35) [fill=black] circle (\vr);
\draw (50,35) [fill=black] circle (\vr);
\draw (20,37.5) node {{\small ${}^{A}$}};
\draw (30,37.5) node {{\small ${}^A$}};
\draw (40,37.5) node {{\small ${}^B$}};
\draw (50,37.5) node {{\small ${}^C$}};
\draw (20,31.5) node {{\small ${}^y$}};
\draw (30,31.5) node {{\small ${}^x$}};
\draw (40,31.5) node {{\small ${}^w$}};
\draw (50,31.5) node {{\small ${}^v$}};
\draw (20,6.5) node {{\small ${}^y$}};
\draw (30,6.5) node {{\small ${}^x$}};
\draw (40,6.5) node {{\small ${}^w$}};
\draw (50,6.5) node {{\small ${}^v$}};
\draw (60,6.5) node {{\small ${}^u$}};
\draw (20,12.5) node {{\small ${}^{B/C}$}};
\draw (30,12.5) node {{\small ${}^A$}};
\draw (40,12.5) node {{\small ${}^A$}};
\draw (50,12.5) node {{\small ${}^B$}};
\draw (60,12.5) node {{\small ${}^C$}};
\draw[densely dashed] (20,10) ellipse (8 and 6);
\draw[densely dashed] (20,35) ellipse (8 and 6);
\draw (15,17.5) node {{\small ${}T$}};
\draw (15,42.5) node {{\small ${}T$}};
\draw (10,28) node {{\small $\cT_1$}};
\draw (10,3) node {{\small $\cT_2$}};
\end{tikzpicture}};
\draw (80,5) node {\begin{tikzpicture}[scale=0.8,style=thick,x=0.1cm,y=0.1cm] \def\vr{2.0pt}\def\vrm{0.5pt}
\draw (20,40)--(70,40);\draw (20,0)--(20,40);
\draw (30,10)--(40,40);\draw (50,10)--(40,40);
\draw (60,0)--(60,40);\draw (70,0)--(70,40);
\draw (20,0) [fill=black] circle (\vr);
\draw (20,10) [fill=black] circle (\vr);
\draw (20,20) [fill=black] circle (\vr);
\draw (20,30) [fill=black] circle (\vr);
\draw (20,40) [fill=black] circle (\vr);
\draw (30,40) [fill=black] circle (\vr);
\draw (30,10) [fill=black] circle (\vr);
\draw (50,10) [fill=black] circle (\vr);
\draw (33.3,20) [fill=black] circle (\vr);
\draw (46.7,20) [fill=black] circle (\vr);
\draw (36.6,30) [fill=black] circle (\vr);
\draw (43.37,30) [fill=black] circle (\vr);
\draw (60,10) [fill=black] circle (\vr);
\draw (60,20) [fill=black] circle (\vr);
\draw (60,30) [fill=black] circle (\vr);
\draw (60,40) [fill=black] circle (\vr);
\draw (60,0) [fill=black] circle (\vr);
\draw (70,10) [fill=black] circle (\vr);
\draw (70,20) [fill=black] circle (\vr);
\draw (70,30) [fill=black] circle (\vr);
\draw (70,40) [fill=black] circle (\vr);
\draw (70,0) [fill=black] circle (\vr);
\draw (40,40) [fill=black] circle (\vr);
\draw (50,40) [fill=black] circle (\vr);
\draw (60,40) [fill=black] circle (\vr);
\draw (70,40) [fill=black] circle (\vr);
\draw (20,42.5) node {{\small ${}^C$}};
\draw (30,42.5) node {{\small ${}^B$}};
\draw (40,42.5) node {{\small ${}^A$}};
\draw (50,42.5) node {{\small ${}^A$}};
\draw (60,42.5) node {{\small ${}^B$}};
\draw (70,42.5) node {{\small ${}^C$}};
\draw (17,-1) node {{\small ${}^C$}};
\draw (17,9) node {{\small ${}^B$}};
\draw (17,19) node {{\small ${}^A$}};
\draw (17,29) node {{\small ${}^A$}};
\draw (57,-1) node {{\small ${}^C$}};
\draw (57,9) node {{\small ${}^B$}};
\draw (57,19) node {{\small ${}^A$}};
\draw (57,29) node {{\small ${}^A$}};
\draw (72.5,-1) node {{\small ${}^C$}};
\draw (72.5,9) node {{\small ${}^B$}};
\draw (72.5,19) node {{\small ${}^A$}};
\draw (72.5,29) node {{\small ${}^A$}};
\draw (33,9) node {{\small ${}^C$}};
\draw (47.3,9) node {{\small ${}^C$}};
\draw (35.8,19) node {{\small ${}^B$}};
\draw (44.1,19) node {{\small ${}^B$}};
\draw (33.8,29) node {{\small ${}^A$}};
\draw (46,29) node {{\small ${}^A$}};
\end{tikzpicture}}; 
\end{tikzpicture}
\caption{Operations $\cT_1$ and $\cT_2$, and a tree belonging to the family $\cT$.}
\label{tree-1}
\end{center}
\end{figure}

If $T \in \cT$, we let $S_A(T)$, $S_B(T)$, and $S_C(T)$ be the sets of vertices of status $A$, $B$,
and $C$, respectively, in $T$. The following observation is immediate from the way in which each tree in the family $\cT$ is constructed.

\begin{observation}[\cite{Haynes-Henning-Hopkins}]
\label{observation1}
If $T \in \cT$, then the following properties hold. \\ [-22pt]
\begin{enumerate}
\item If $v\in S_A(T)$, then $|N_T(v)\cap (S_B(T)\cup S_C(T))|=1$ and $N_T(v)\cap S_A(T)\not=\emptyset$; \1
\item If $v \in S_B(T)$ $($$v\in S_C(T)$, resp.$)$, then $|N_T(v)\cap S_C(T)|=1$ $($$|N_T(v)\cap S_B(T)|=1$, resp.$)$, and $N_T(v) \setminus (S_B(T)\cup S_C(T)) \subseteq S_A(T)$; \1
\item $L(T)\subseteq S_C(T)$ and $S(T)\subseteq S_B(T)$; \1
\item $\{v\in V(T)\colon \max\{d_T(v,S_B(T)), d_T(v,S_C(T))\}=2\}
 = S_A(T)$; \1
\item $|S_B(T)|=|S_C(T)|$.
\end{enumerate}
\end{observation}

Haynes, Henning and Hopkins~\cite{Haynes-Henning-Hopkins} gave the following characterization of trees in Class~$3$.

\begin{thm}[\cite{Haynes-Henning-Hopkins}]
\label{theorem of Haynes-Henning-Hopkins}
A tree is in Class~$3$ if and only if it is the underlying tree of a labeled tree that belongs to the family~$\cT$.
\end{thm}

\section{Two new characterizations}

We are interested in structural properties of trees belonging to Class~$3$. We first present another constructive characterization of labeled trees belonging to the family $\cT$ that is a modification of the characterization given in~\cite{Haynes-Henning-Hopkins}.

\begin{defn}
\label{defn:cO}
Let $\cO$ be the family of labeled trees that:
\begin{itemize}
\item[$(1)$] contains a path $P_6$ in which the two central vertices have status A, and all other vertices have status B; and
\item[$(2)$] is closed under the two operations $\cO_1$ and $\cO_2$, which extend the labeled tree $T$ belonging to $\cO$ by adding a labeled $P_6$, and then:
\end{itemize}
\begin{itemize}
\item {\bf Operation} $\cO_1$. Identifying one B-B-A-path of $T$ with one B-B-A-path of $P_6$;
\item {\bf Operation} $\cO_2$. Identifying one B-B-edge of $T$ with one B-B-edge of $P_6$.
\end{itemize}
\end{defn}

The two operations $\cO_1$ and $\cO_2$, and an example of a labeled tree belonging to the family $\cO$ are given in Fig.~\ref{tree-2}.

\begin{figure}[htb] \begin{center} 
\begin{tikzpicture}[scale=0.8,style=thick,x=0.1cm,y=0.1cm]
\def\vr{2.5pt}
\draw (0,0) node {\begin{tikzpicture}[scale=1,style=thick,x=0.1cm,y=0.1cm] \def\vr{2.0pt}
\draw (20,10)--(70,10); \draw (20,40)--(70,40);
\draw (20,10) [fill=black] circle (\vr);
\draw (30,10) [fill=black] circle (\vr);
\draw (40,10) [fill=black] circle (\vr);
\draw (50,10) [fill=black] circle (\vr);
\draw (60,10) [fill=black] circle (\vr);
\draw (70,10) [fill=black] circle (\vr);
\draw (20,40) [fill=black] circle (\vr);
\draw (30,40) [fill=black] circle (\vr);
\draw (40,40) [fill=black] circle (\vr);
\draw (50,40) [fill=black] circle (\vr);
\draw (60,40) [fill=black] circle (\vr);
\draw (70,40) [fill=black] circle (\vr);
\draw (20,42.5) node {{\small ${}^B$}};
\draw (30,42.5) node {{\small ${}^B$}};
\draw (40,42.5) node {{\small ${}^A$}};
\draw (50,42.5) node {{\small ${}^A$}};
\draw (60,42.5) node {{\small ${}^B$}};
\draw (70,42.5) node {{\small ${}^B$}};
\draw (20,12.5) node {{\small ${}^B$}};
\draw (30,12.5) node {{\small ${}^B$}};
\draw (40,12.5) node {{\small ${}^A$}};
\draw (50,12.5) node {{\small ${}^A$}};
\draw (60,12.5) node {{\small ${}^B$}};
\draw (70,12.5) node {{\small ${}^B$}};
\draw[densely dashed] (25,10) ellipse (12 and 8);
\draw[densely dashed] (30,40) ellipse (17 and 10);
\draw (15,19.5) node {{\small ${}T$}};
\draw (15,47.5) node {{\small ${}T$}};
\draw (10,33) node {{\small $\cO_1$}};
\draw (10,0) node {{\small $\cO_2$}};
\end{tikzpicture}};
\draw (90,0) node {\begin{tikzpicture}[scale=0.8,style=thick,x=0.1cm,y=0.1cm] \def\vr{2.0pt}\def\vrm{0.5pt}
\draw (20,40)--(70,40);\draw (20,0)--(20,40);
\draw (30,10)--(40,40);\draw (50,10)--(40,40);
\draw (60,0)--(60,40);\draw (70,0)--(70,40);
\draw (20,0) [fill=black] circle (\vr);
\draw (20,10) [fill=black] circle (\vr);
\draw (20,20) [fill=black] circle (\vr);
\draw (20,30) [fill=black] circle (\vr);
\draw (20,40) [fill=black] circle (\vr);
\draw (30,40) [fill=black] circle (\vr);
\draw (30,10) [fill=black] circle (\vr);
\draw (50,10) [fill=black] circle (\vr);
\draw (33.3,20) [fill=black] circle (\vr);
\draw (46.7,20) [fill=black] circle (\vr);
\draw (36.6,30) [fill=black] circle (\vr);
\draw (43.37,30) [fill=black] circle (\vr);
\draw (60,10) [fill=black] circle (\vr);
\draw (60,20) [fill=black] circle (\vr);
\draw (60,30) [fill=black] circle (\vr);
\draw (60,40) [fill=black] circle (\vr);
\draw (60,0) [fill=black] circle (\vr);
\draw (70,10) [fill=black] circle (\vr);
\draw (70,20) [fill=black] circle (\vr);
\draw (70,30) [fill=black] circle (\vr);
\draw (70,40) [fill=black] circle (\vr);
\draw (70,0) [fill=black] circle (\vr);
\draw (40,40) [fill=black] circle (\vr);
\draw (50,40) [fill=black] circle (\vr);
\draw (60,40) [fill=black] circle (\vr);
\draw (70,40) [fill=black] circle (\vr);
\draw (20,42.5) node {{\small ${}^B$}};
\draw (30,42.5) node {{\small ${}^B$}};
\draw (40,42.5) node {{\small ${}^A$}};
\draw (50,42.5) node {{\small ${}^A$}};
\draw (60,42.5) node {{\small ${}^B$}};
\draw (70,42.5) node {{\small ${}^B$}};
\draw (17,-1) node {{\small ${}^B$}};
\draw (17,9) node {{\small ${}^B$}};
\draw (17,19) node {{\small ${}^A$}};
\draw (17,29) node {{\small ${}^A$}};
\draw (57,-1) node {{\small ${}^B$}};
\draw (57,9) node {{\small ${}^B$}};
\draw (57,19) node {{\small ${}^A$}};
\draw (57,29) node {{\small ${}^A$}};
\draw (72.5,-1) node {{\small ${}^B$}};
\draw (72.5,9) node {{\small ${}^B$}};
\draw (72.5,19) node {{\small ${}^A$}};
\draw (72.5,29) node {{\small ${}^A$}};
\draw (32.5,9) node {{\small ${}^B$}};
\draw (47.3,9) node {{\small ${}^B$}};
\draw (35.8,19) node {{\small ${}^B$}};
\draw (44.1,19) node {{\small ${}^B$}};
\draw (33.8,29) node {{\small ${}^A$}};
\draw (46,29) node {{\small ${}^A$}};
\end{tikzpicture}};
\end{tikzpicture}
\caption{Operations $\cO_1$ and $\cO_2$, and a tree belonging to the family $\cO$.} \label{tree-2}
\end{center}
\end{figure}

If $T\in \cO$, we let $S_A(T)$ and $S_B(T)$ be the sets of vertices of status $A$ and $B$, respectively, in~$T$. The following observation follows immediately from the way in which each tree in the family~$\cO$ is constructed.

\begin{observation}
\label{observation2} If $T \in \cO$, then the following two properties hold. \\ [-22pt]
\begin{enumerate}
\item $|N_T(v)\cap S_B(T)|=1$ for each $v\in V(T)$; \1
\item $L(T)\cup S(T) \subseteq S_B(T)$.
\end{enumerate}
\end{observation}

We are now in a position to prove the following characterization of trees that belong to Class~$3$.

\begin{thm}
\label{new-characterization}
If $T$ is a tree of order at least~$6$, then the following statements are equivalent:
\\ [-22pt]
\begin{enumerate}
\item[{\rm (1)}] $T$ is in Class~$3$; \1
\item[{\rm (2)}]  $T \in \cT$; \1
\item[{\rm (3)}]  $T \in \cO$; \1
\item[{\rm (4)}]  There is a uniquely determined subset $F$ of $E(T)$ such that:
\begin{itemize}
\item[{\rm (a)}]  each pendant edge of $T$ belongs to $F$;
\item[{\rm (b)}]  $d_T(e,F \setminus \{e\})=3$ for each $e\in F$;
\item[{\rm (c)}]  if $e$ and $f$ are distinct elements of $F$, then there is a unique sequence $(e_0, e_1, \ldots, e_k)$ of elements of $F$ such that $e_0=e$, $e_k=f$, and $d_T(e_{i-1},e_i)=3$ for $i \in [k]$.
\end{itemize}
\end{enumerate}
\end{thm}

\begin{figure}[htb]
\begin{center}
\begin{tikzpicture}[scale=0.7,style=thick,x=0.1cm,y=0.1cm]
\def\vr{2.5pt}
\draw (0,0)--(90,0); \draw (70,10)--(100,10);
\draw (30,0)--(30,30); \draw (70,0)--(70,30);
\draw[line width=2.5pt] (0,0)--(10,0);
\draw[line width=2.5pt] (40,0)--(50,0);
\draw[line width=2.5pt] (80,0)--(90,0);
\draw[line width=2.5pt] (90,10)--(100,10);
\draw[line width=2.5pt] (30,20)--(30,30);
\draw[line width=2.5pt] (70,20)--(70,30);
\draw (0,0) [fill=black] circle (\vr);
\draw (10,0) [fill=black] circle (\vr);
\draw (20,0) [fill=black] circle (\vr);
\draw (30,0) [fill=black] circle (\vr);
\draw (40,0) [fill=black] circle (\vr);
\draw (50,0) [fill=black] circle (\vr);
\draw (60,0) [fill=black] circle (\vr);
\draw (70,0) [fill=black] circle (\vr);
\draw (80,0) [fill=black] circle (\vr);
\draw (90,0) [fill=black] circle (\vr);
\draw (30,10) [fill=black] circle (\vr);
\draw (30,20) [fill=black] circle (\vr);
\draw (30,30) [fill=black] circle (\vr);
\draw (70,10) [fill=black] circle (\vr);
\draw (70,20) [fill=black] circle (\vr);
\draw (70,30) [fill=black] circle (\vr);
\draw (80,10) [fill=black] circle (\vr);
\draw (90,10) [fill=black] circle (\vr);
\draw (100,10) [fill=black] circle (\vr);
\draw (5,2) node {${}^{e_0=e}$};
\draw (45,2) node {${}^{e_1}$};
\draw (85,2) node {${}^{e_2}$};
\draw (73.5,24) node {${}^{e_3}$};
\draw (95,12.5) node {${}^{e_4=f}$};
\end{tikzpicture}
\caption{An example to the statement (d) of Theorem~\ref{new-characterization}.}
\label{an-example}
\end{center}
\end{figure}
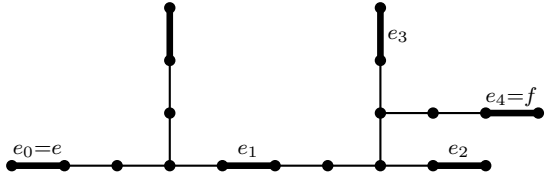

For example, if $T$ is the tree illustrated in Fig.~\ref{an-example}, then the subset $F$ of broad edges of $T$ indicated in Fig.~\ref{an-example} satisfies property~(4) in the statement of Theorem~\ref{new-characterization}. Moreover if $e = e_0$ and $f = e_4$, then the sequence $(e_0, e_1,e_2,e_3,e_4)$ satisfies property~(4c) in the statement of Theorem~\ref{new-characterization}.

\begin{proof} The equivalence of the statements $(1)$ and $(2)$
has been proved in \cite{Haynes-Henning-Hopkins}, see Theorem~\ref{theorem of Haynes-Henning-Hopkins}. The proof of the equivalence of the statements (2) and (3) is straightforward by Observations~\ref{observation1} and \ref{observation2}, and we omit the details. We shall prove the equivalence of the statements $(3)$ and $(4)$.

Assume first that $T \in \cO$. Thus the tree $T$ can be obtained from a sequence $T_1,\ldots,T_p$ of trees, where $T_1=P_6$ and $T=T_p$, and, if $p \ge 2$, $T_{i+1}$ can be obtained from $T_i$ by operation $\cO_1$ or $\cO_2$ for $i \in [p-1]$. Let $F$ be the set of all B-B edges in $T$. By induction on the number~$p$ we shall prove that $F$ has the desired properties.

If $p=1$, then $T=P_6$, and the set $F$ (consisting of the two pendant edges of $P_6$) has the desired properties. This establishes the base case. Assume, then, that the result holds for all trees that can be constructed from a~sequence of fewer than $p$ trees, where $p\ge 2$. Let $T \in \cO$ be obtained from a~sequence $T_1,T_2,\ldots,T_p$ of $p$ trees. By our inductive hypothesis, the set $F'$ of the $B$-$B$ edges of $T_{p-1}$ has the desired properties in $T_{p-1}$. We now consider two possibilities depending on whether $T$ is obtained from $T_{p-1}$ by operation $\cO_1$ or $\cO_2$.

Assume first that $T$ is obtained from $T_{p-1}$ by operation $\cO_1$. Suppose, without loss of generality, that $T$ is obtained from $T_{p-1}$ by identifying a B-B-A path $xyz$ in $T_{p-1}$ with a B-B-A path in the labeled $P_6$, say with $xyz$ in $P_6 \colon xyzabc$, where $z$ and $a$ have status A, while $x, y, b, c$ have status B. Now from the properties of $F'$ in $T_{p-1}$, we infer that the set $F' \cup \{bc\}$ (of all B-B edges in $T$) has  properties (4a)--(4c) in $T$. Similarly, if $T$ is obtained from $T_{p-1}$ by operation $\cO_2$, then again the set of all B-B edges of $T$ has properties (4a)--(4c) in $T$.

Assume now that $T$ is a tree of order at least~$6$ in which there is a unique nonempty subset $F$ of $E(T)$ having properties $({\rm 4a})$--$({\rm 4c})$. By induction on the order $n \ge 6$ of $T$, we shall prove that $T \in \cO$. The only tree of order~$6$ having the desired properties is the path $P_6$, and so in this case $T \in \cO$. This establishes the base step. Let $n>6$, and assume that if $T'$ is a tree of order $n'$ with $6 \le n'<n$ that has a unique subset of edges having properties $({\rm 4a})$--$({\rm 4c})$, then $T' \in \cO$. Let $T$ be a tree of order $n$ with a unique subset $F$ of $E(T)$ having properties $({\rm 4a})$--$({\rm 4c})$. Let $P \colon v_0v_1 \ldots  v_k$ be a longest path in $T$. For convenience, we root $T$ at the leaf~$v_k$.
From the properties (4a) and (4b) of edges belonging to $F$ and by the maximality of the path $P$, it follows that $k \ge 5$ and $d_T(v_1)=d_T(v_2)=2$. We consider two cases depending on the degree of $v_3$ in $T$.

\smallskip
\emph{Case~1. $d_T(v_3) = 2$.}  In this case, we let $T'$ denote the subtree $T - \{v_0, v_1, v_2, v_3\}$ of $T$. Since $v_0v_1$  is a pendant edge in $T$, it follows from the properties (4a) and (4b) that there is exactly one vertex $v_4' \in N_T(v_4) \setminus \{v_3\}$ such that $v_4v_4'\in F$. Now, from the properties $({\rm 4a})$--$({\rm 4c})$ of $F$ in $T$ it follows that $F'=F \setminus \{v_0v_1\}$ has the properties $({\rm 4a})$--$({\rm 4c})$ in $T'$. Applying the inductive hypothesis we infer that $T' \in \cO$. Thus, $T$ can be obtained from $T'$ by operation $\cO_2$ (identifying the B-B edge $v_4v_4'$ of $T'$ with the B-B edge $v_4v_4'$ of the labeled $P_6 \colon v_0v_1v_2v_3v_4v_4'$), and so, $T\in \cO$.

\smallskip
\emph{Case~2. $d_T(v_3) \ge 3$.} In this case, we let $T'$ denote the subtree $T-\{v_0, v_1, v_2\}$ of $T$. Since $v_0v_1$  is a pendant edge in $T$, it follows from the properties (4a) and (4b) that there are vertices $x \in N_T(v_3) \setminus \{v_2\}$ and $x' \in N_T(x) \setminus \{v_3\}$ such that $xx' \in F$. As in Case 1, from the properties $({\rm 4a})$--$({\rm 4c})$ of $F$ in $T$ it follows that $F'=F-\{v_0v_1\}$ has the properties $({\rm 4a})$--$({\rm 4c})$ in $T'$. Hence, an application of the inductive hypothesis implies that $T'\in \cO$.
In this case, $T$ can be obtained from $T'$ by applying operation $\cO_1$ (identifying the B-B-A path $x'xv_3$ of $T'$ with the B-B-A path $x'xv_3$ of the labeled $P_6 \colon v_0v_1v_2v_3xx'$). This proves that $T\in \cO$, and completes the proof of the equivalence of the statements (3) and (4).
\end{proof}

\smallskip
As an immediate consequence of Theorem~\ref{new-characterization}, we characterize the paths that are in Class~$3$.

\begin{cor} \label{wniosek-path-in-Class-3}
A path $P_n$ is in Class~$3$ if and only if $n \equiv 2 \, (\modo \, 4)$.
\end{cor}

\begin{figure}[htb] \begin{center} 
\begin{tikzpicture}[scale=0.8,style=thick,x=0.1cm,y=0.1cm] \def\vr{2.6pt}
\draw (-20,0) node {\begin{tikzpicture}[scale=0.8,style=thick,x=0.1cm,y=0.1cm] \def\vr{2.0pt}
\draw (-20,28.3)--(0,20)--(-20,11.7);
\draw (0,20)--(0,35);
\draw (0,20) [fill=black] circle (\vr);
\draw (0,35) [fill=black] circle (\vr);
\draw[fill=white] (-25,20) ellipse (10 and 10);
\draw[fill=white] (-20,20) ellipse (5 and 8.3);
\draw (-15,31) node {\small $T$};
\draw (0,16.5) node {\small ${}^v$};
\draw (-20,15.5) node {{\small ${}^A$}};
\end{tikzpicture}}; 
\draw (13.8,0) node {\begin{tikzpicture}[scale=0.8,style=thick,x=0.1cm,y=0.1cm] \def\vr{2.0pt}
\draw (20,28.3)--(0,20)--(20,11.7);
\draw (0,20) [fill=black] circle (\vr);
\draw (0,35) [fill=black] circle (\vr);
\draw[fill=white] (25,20) ellipse (10 and 10);
\draw[fill=white] (20,20) ellipse (5 and 8.3);
\draw (20,15.5) node {{\small ${}^B$}};
\end{tikzpicture}}; 
\draw (65,0) node
{\begin{tikzpicture}[scale=0.8,style=thick,x=0.1cm,y=0.1cm] \def\vr{2.0pt}
\draw (-20,28.3)--(0,20)--(-20,11.7);
\draw (0,20)--(0,35);
\draw (0,20) [fill=black] circle (\vr);
\draw (0,35) [fill=black] circle (\vr);
\draw[fill=white] (-25,20) ellipse (10 and 10);
\draw[fill=white] (-20,20) ellipse (5 and 8.3);
\draw (-15,31) node {\small $T_A$};
\draw (0,16.5) node {\small ${}^v$};
\draw (-20,15.5) node {{\small ${}^A$}};
\end{tikzpicture}};
\draw (110,0) node
{\begin{tikzpicture}[scale=0.8,style=thick,x=0.1cm,y=0.1cm] \def\vr{2.0pt}
\draw (20,28.3)--(0,20)--(20,11.7);
\draw (0,20)--(0,35);
\draw (0,20) [fill=black] circle (\vr);
\draw (0,35) [fill=black] circle (\vr);
\draw[fill=white] (25,20) ellipse (10 and 10);
\draw[fill=white] (20,20) ellipse (5 and 8.3);
\draw (15,31) node {\small $T_B$};
\draw (0,16.5) node {\small ${}^v$};
\draw (20,15.5) node {{\small ${}^B$}};
\end{tikzpicture}};  
\draw (-2.5,-20) node {{\small (a) $T$}};
\draw (70,-20) node {{\small (b) $T_A$}};
\draw (110,-20) node {{\small (c) $T_B$}};
\end{tikzpicture}
\caption{A tree $T$ and the components $T_A$ and $T_B$} \label{tree-operacja11111}
\end{center}
\end{figure}
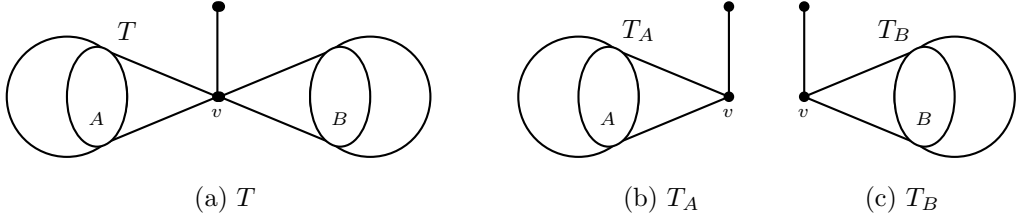

\newpage
Let $v$ be a weak support vertex of degree at least~$3$ in a tree $T$. If sets $A$ and $B$ form a partition of the set $N_T(v) \setminus L(T)$, then by $T_A$ ($T_B$, resp.) we denote the component of $T-B$ ($T-A$, resp.) that contains the vertex $v$, as illustrated in Fig.~\ref{tree-operacja11111}.

\begin{lem}
\label{wniosek-o-TA-i-TB}
Let $v$ be a weak support vertex of degree at least~$3$ in a tree $T$, and let sets $A$ and $B$ form a partition of the set $N_T(v) \setminus L(T)$. Then the tree $T$ is in Class~$3$ if and only if both subtrees $T_A$ and $T_B$ are in Class~$3$.
\end{lem}
\begin{proof}
The proof follows readily from the fact that a subset $F$ of $E(T)$ has properties (4a)--(4c) of Theorem~\ref{new-characterization} in $T$ if and only if each of the sets $F \cap E(T_A)$ and $F \cap E(T_B)$ has properties (4a)--(4c) of Theorem~\ref{new-characterization} in $T_A$ and $T_B$, respectively.
\end{proof}

A \emph{caterpillar} is a tree of order at least~$3$ with the property that the removal of its leaves results in a path, called the \emph{spine} of the caterpillar. The code $C$ of a caterpillar $T$ with spine $v_0v_1 \ldots v_s$ is the sequence of nonnegative integers $(t_0, t_1, \ldots, t_s)$, where $t_i$ is the number of leaves adjacent to $v_i$ in $T$. We say that two leaves of a caterpillar $T$ are \emph{consecutive} if no inner vertex of the path joining their neighbors is a support vertex in $T$. Haynes et al.~\cite{Haynes-Hedetniemi-van der Merwe} characterized caterpillars in Class~$3$.

\begin{thm}[\cite{Haynes-Hedetniemi-van der Merwe}]
\label{thm:caterpillar}
A caterpillar $T$ with code $C=(t_0, \ldots, t_s)$ is in Class~$3$ if and only if $t_i \in \{0,1\}$ for $i \in [s]_0$, and any two consecutive nonzero entries in $C$ are at distance $3 \, (\modo \, 4)$.
\end{thm}

Let $T$ be a caterpillar with spine $v_0v_1 \ldots v_s$ and with code $C=(t_0, \ldots, t_s)$ that is in Class~$3$. If $t_i$ and $t_j$ are two consecutive nonzero entries in $C$, then by Theorem~\ref{thm:caterpillar} we have $t_i = t_j = 1$ and the vertices $v_i$ and $v_j$ on the spine are at distance $3 \, (\modo \, 4)$, that is, $d_T(v_i,v_j) \equiv \, 3 \, (\modo \, 4)$. Let $l_i$ and $l_j$ be the leaf neighbors of the vertices $v_i$ and $v_j$, respectively. Then, $l_i$ and $l_j$ are consecutive leaves in $T$ and $d_T(l_i,l_j) \equiv \, 1 \, (\modo \, 4)$. Conversely, if $T$ is a caterpillar with code $C=(t_0, \ldots, t_s)$ where $t_i \in \{0,1\}$ and where any two consecutive leaves of $T$ are at distance $1 \, (\modo \, 4)$, then the neighbors of these leaves (that belong to the spine of $T$) are at distance $3 \, (\modo \, 4)$, and so by Theorem~\ref{thm:caterpillar}, the caterpillar $T$ is in Class~$3$. This yields the following equivalent statement of Theorem~\ref{thm:caterpillar}.

\begin{thm}[\cite{Haynes-Hedetniemi-van der Merwe}]
\label{thm:caterpillar2}
A caterpillar $T$ of order at least~$6$ is in Class~$3$ if and only if every two consecutive leaves of $T$ are at distance $1 \, (\modo \, 4)$.
\end{thm}

We present next a short proof of Theorem~\ref{thm:caterpillar2}.

\noindent
\begin{proof}
Let $T$ be a caterpillar of order~$n \ge 6$. We proceed by induction on the number $l \ge 2$ of leaves in $T$ to show that $T$ is in Class~$3$ if and only if every two consecutive leaves of $T$ are at distance $1 \, (\modo \, 4)$. If $l=2$, then $T$ is a path~$P_n$, and so by Corollary~\ref{wniosek-path-in-Class-3}, $n \equiv 2 \, (\modo \, 4)$, and so the two leaves of $T$ are at distance $1 \, (\modo \, 4)$. This establishes the base case. Let $l \ge 3$ and assume that every caterpillar $T'$ of order at least~$6$ with $l'$ leaves is in Class~$3$ if and only if every two consecutive leaves of $T'$ are at distance $1 \, (\modo \, 4)$. Let $T$ be a caterpillar of order~$n \ge 6$ with $l$ leaves.

If $T$ has a strong support vertex $v$ and if $v_1$ and $v_2$ are two leaf neighbors of $v$,  then subdividing two edges $vv_1$ and $vv_2$ increases the total domination number, implying that $T$ is in Class~$1$ or~$2$, a contradiction. Hence, every support vertex of $T$ is a weak support vertex. Since $T$ has at least three leaves, the caterpillar $T$ has a leaf, say $d$, whose (unique) neighbor, say $c$, is of degree~$3$. Let $a$ and $b$ be the two neighbors of $c$ on the spine of $T$. Let $T_a$ ($T_b$, resp.) be the component of $T-b$ ($T-a$, resp.) that contains $a$ ($b$, resp.). We note that both $T_a$ and $T_b$ are caterpillars. By Corollary~\ref{wniosek-o-TA-i-TB}, $T$ is in Class~$3$ if and only if both $T_a$ and $T_b$ are in Class~$3$. By the inductive hypothesis, $T_a$ and $T_b$ are in Class~$3$ if and only if every two consecutive leaves in $T_a$ are at distance $1 \, (\modo \, 4)$ and every two consecutive leaves in $T_b$ are at distance $1 \, (\modo \, 4)$. From these properties of the caterpillars $T_a$ and $T_b$, we infer that every two consecutive leaves of $T$ are at distance $1 \, (\modo \, 4)$. This proves Theorem~\ref{thm:caterpillar2}.
\end{proof}

\section{Another characterization}

We show in this section that trees belonging to the family $\cO$ are $2$-subdivisions of trees or can be obtained from $2$-subdivisions of trees, where a $2$-subdivision of a graph is defined as follows. Recall that the \emph{corona} $H \circ K_1$ of a graph $H$ is the graph obtained from $H$ by adding for each vertex $v \in V(H)$ a new vertex $v'$ and the edge $vv'$.

\begin{defn}
\label{defn:2subdivision}
Let $G$ be a connected graph of order at least~$2$, and let ${\cal P} = \{{\cal P}(v) \colon v \in V(G)\}$ be a family in which ${\cal P}(v)$ is a partition of the set $N_G(v)$ for each $v \in V(G)$. The $2$-subdivision of $G$ with respect to ${\cal P}$ is the graph $G({\cal P})$ with vertex set
\[
V(G({\cal P})) = V(G)\cup \left(V(G)\times \{1\}\right) \cup \bigcup_{v\in V(G)}\left(\{v\}\times {\cal P}(v)\right)
\]
and edge set $E(G({\cal P})) = E_1 \cup E_2 \cup E_3$ where
\[
\begin{array}{lcl}
E_1 & = & \displaystyle{ \{v(v,1)\colon v\in V(G)\} }, \2 \\
E_2 & = & \displaystyle{ \bigcup_{v\in V(G)} \{v(v,A)\colon A\in {\cal P}(v)\} }, \2 \\
E_3 & = & \displaystyle{ \bigcup_{uv\in E(G)}\{(u,A)(v,B)\colon A\in {\cal P}(u), B\in {\cal P}(v), u\in B, v\in A\}  }. \\
\end{array}
\]
More intuitively, $G({\cal P})$ is the graph obtained from the corona $G\circ K_1$ by inserting two new vertices into each inner edge of $G\circ K_1$, and then identifying newly inserted vertices according to the partition ${\cal P}(v)$ of $N_G(v)$, that is, if $A=\{w_1,\ldots,w_k\}\in {\cal P}(v)$, then we contract all neighbors of $v$ on the edges $vw_1,\ldots, vw_k$ into a single vertex $(v,A)$ and replace all multiple edges in the resulting graph by single edges, for each $v\in V(G)$ and $A\in {\cal P}(v)$. 
\end{defn}

\begin{figure}[htb] \begin{center} 
\begin{tikzpicture}[scale=0.7,style=thick,x=0.1cm,y=0.1cm] \def\vr{3pt}
\draw (-15,-5.4) node {\begin{tikzpicture}[scale=0.7,style=thick,x=0.1cm,y=0.1cm] \def\vr{3pt}
\draw (60,10)--(0,10)--(70,35)--(10,35);
\draw (70,35)--(20,60)--(80,60);
\draw (0,10) [fill=black] circle (\vr);
\draw (60,10) [fill=black] circle (\vr);
\draw (10,35) [fill=black] circle (\vr);
\draw (70,35) [fill=black] circle (\vr);
\draw (20,60) [fill=black] circle (\vr);
\draw (80,60) [fill=black] circle (\vr);
\draw (20,56) node {\small ${}^a$};
\draw (10,31) node {\small ${}^b$};
\draw (0,6) node {\small ${}^c$};
\draw (80,56) node {\small ${}^d$};
\draw (70,31) node {\small ${}^e$};
\draw (60,5.6) node {\small ${}^f$};
\draw (0,60) node {\small $G$};
\end{tikzpicture}};
\draw (90,0) node {\begin{tikzpicture}[scale=0.7,style=thick,x=0.1cm,y=0.1cm] \def\vr{3pt}
\draw (0,10)--(0,22.5);\draw (0,22.5) [fill=black] circle (\vr);
\draw (10,35)--(10,47.5);\draw (10,47.5) [fill=black] circle (\vr);
\draw (40,60)--(20,60)--(20,72.5);\draw (20,72.5) [fill=black] circle (\vr);
\draw (60,10)--(60,22.5);\draw (60,22.5) [fill=black] circle (\vr);
\draw (70,35)--(70,47.5);\draw (70,47.5) [fill=black] circle (\vr);
\draw (80,60)--(80,72.5);\draw (80,72.) [fill=black] circle (\vr);
\draw (60,10)--(0,10); \draw (20,10)--(70,35)--(10,35);
\draw (50,35)--(40,60)--(80,60);
\draw (0,10) [fill=black] circle (\vr);
\draw (20,10) [fill=black] circle (\vr);
\draw (40,10) [fill=black] circle (\vr);
\draw (10,35) [fill=black] circle (\vr);
\draw (30,35) [fill=black] circle (\vr);
\draw (50,35) [fill=black] circle (\vr);
\draw (20,60) [fill=black] circle (\vr);
\draw (40,60) [fill=black] circle (\vr);
\draw (60,60) [fill=black] circle (\vr);
\draw (60,10) [fill=black] circle (\vr);
\draw (10,35) [fill=black] circle (\vr);
\draw (70,35) [fill=black] circle (\vr);
\draw (20,60) [fill=black] circle (\vr);
\draw (80,60) [fill=black] circle (\vr);
\draw (45,22.5) [fill=black] circle (\vr);
\draw (20,56) node {\small ${}^a$};
\draw (13,72.5) node {\small ${}^{(a,1)}$};
\draw (87,72.5) node {\small ${}^{(d,1)}$};
\draw (10,31) node {\small ${}^b$};
\draw (3,47.5) node {\small ${}^{(b,1)}$};
\draw (77,47.5) node {\small ${}^{(e,1)}$};
\draw (0,6) node {\small ${}^c$};
\draw (-7,22.5) node {\small ${}^{(c,1)}$};
\draw (67,22.5) node {\small ${}^{(f,1)}$};
\draw (80,56) node {\small ${}^d$};
\draw (70,31) node {\small ${}^e$};
\draw (60,5.6) node {\small ${}^f$};
\draw (0,60) node {\small $G({\cal P})$};
\draw (40,64) node {\small ${}^{(a,\{d,e\})}$};
\draw (60,64) node {\small ${}^{(d,\{a\})}$};
\draw (30,39) node {\small ${}^{(b,\{e\})}$};
\draw (49,30.5) node {\small ${}^{(e,\{a,b\})}$};
\draw (36.5,23.5) node {\small ${}^{(e,\{c\})}$};
\draw (40.5,5) node {\small ${}^{(f,\{c\})}$};
\draw (20,5) node {\small ${}^{(c,\{e,f\})}$};
\end{tikzpicture}};
\draw (-20,-45) node {{\small (a) $G$}};
\draw (90,-45) node {{\small (b) $G({\cal P})$}};
\end{tikzpicture}
\caption{The trees $G$ and $G({\cal P})$.} \label{rys-do-podrodziny}
\end{center}
\end{figure}
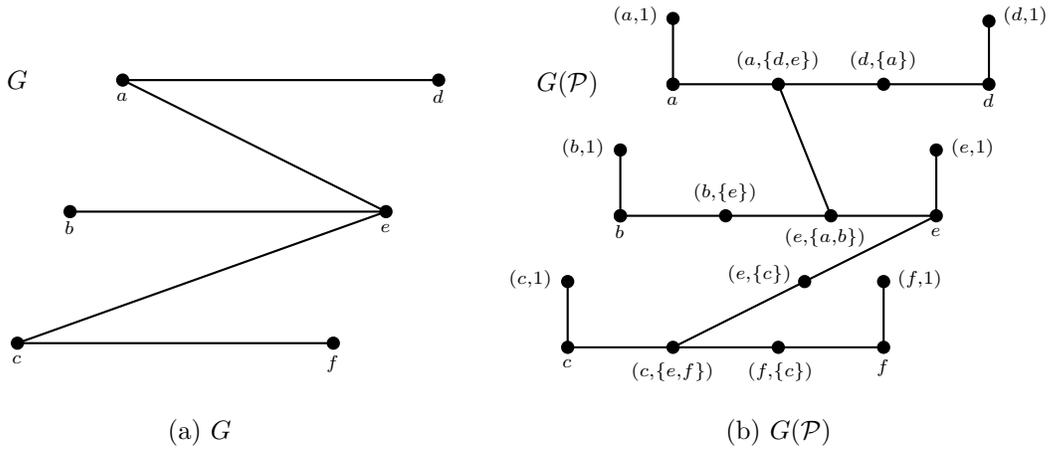

\begin{ex}
If $G$ is the tree shown in Fig.~\ref{rys-do-podrodziny}(a), and ${\cal P}= \{{\cal P}(a), \ldots, {\cal P}(f)\}$ is a~family of partitions of the sets $N_G(a), \ldots, N_G(f)$, respectively, where ${\cal P}(a)=  \{\{d, e\}\}$, ${\cal P}(b)= \{\{e\}\}$, ${\cal P}(c) =\{\{e, f\}\}$, ${\cal P}(d)=\{\{a\}\}$, ${\cal P}(e)=\{\{a, b\},\{c\}\}$, and ${\cal P}(f)=\{\{c\}\}$, then $G({\cal P})$ is the tree shown in Fig.~\ref{rys-do-podrodziny}(b).
\end{ex}

Let $\cO'$ be the family of all trees $T({\cal P})$, where $T$ is any tree of order at least~$2$, and ${\cal P} = \{{\cal P}(v)\colon v\in V(T)\}$ is a~family in which ${\cal P}(v)$ is a~partition of the set $N_T(v)$ for each $v\in V(T)$. If $T({\cal P}) \in \cO'$, then we observe that the set $F = \{v(v,1) \colon v\in V(T)\}$  of all pendant edges of $T({\cal P})$ has properties (4a)--(4c) of Theorem~\ref{new-characterization}. Therefore,  $T({\cal P})\in \cO$, and so $\cO' \subseteq \cO$. The path $P_{10}$ proves that $\cO'$ is a proper subfamily of the family $\cO$. We state these observations formally as follows.

\begin{observation}
The family $\cO'$ is a proper subfamily of the family $\cO$.
\end{observation}

In the next theorem we characterize trees belonging to the family $\cO'$ in terms of $2$-packings. A set $S$ of vertices of a graph $G$ is a $2$-\emph{packing} in $G$ if the vertices in $S$ are pairwise at distance at least~$3$ in $G$, that is, $N_G[u] \cap N_G[v] = \emptyset$ for every pair of distinct vertices $u,v \in S$. The $2$-\emph{packing number} of $G$, denoted by $\rho(G)$, is the maximum cardinality of a $2$-packing in $G$.

\begin{thm}
\label{obserwacja-O-prim-$2$-packing}
If $T$ is a tree of order at least~$6$, then $T$ is in the family $\cO'$ if and only if the set of weak support vertices of $T$ is a maximum $2$-packing in $T$.
\end{thm}
\begin{proof}
Assume first that $T\in \cO'$, say $T=R({\cal P})$ for some tree $R$ and some family ${\cal P}$ of partitions of the sets $N_R(v)$, $v\in V(R)$. From the definition of $R({\cal P})$ it follows that $V(R)$ is the set of weak supports in $R({\cal P})$. In addition, the set $V(R)$ is also a $2$-packing in $R({\cal P})$ (as $d_{R({\cal P})}(u,v)\ge 3$ for each pair of vertices $u, v\in V(R)$), and therefore $|V(R)|\le \rho(R({\cal P}))$. On the other hand since $\{N_{R({\cal P})}[v]\colon v\in V(R)\}$ is a partition of the set $V(R({\cal P}))$ and every $2$-packing has at most one vertex in $N_R[v]$ for every $v\in V(R)$ (as $d_{R({\cal P})}(x,y)\le 2$ for each pair of vertices $x, y\in N_{R({\cal P})}[v]$) it follows that every $2$-packing in $R({\cal P})$ has at most $|V(R)|$ vertices, and therefore $\rho(R({\cal P}))\le |V(R)|$. This implies that $V(R)$ is a maximum $2$-packing in $R({\cal P})$.

Assume now that $T$ is a tree of order at least~$6$ in which the set $S'(T)$ of weak support vertices is a maximum $2$-packing. We shall prove that $T$ is in the family $\cO'$. We first prove three claims.

\smallskip
\begin{claim}
\label{claim1}
The tree $T$ has no strong support vertex.
\end{claim}
\noindent \emph{Proof of Claim~\ref{claim1}}.
Suppose, to the contrary, that $v$ is a strong support vertex in $T$. Let $v'$ be any leaf adjacent to $v$ in $T$. Certainly, neither $v$ nor $v'$ is in $S'(T)$. Since $S'(T)$ is a maximum $2$-packing in $T$, the set $S'(T)\cup \{v'\}$ is not a $2$-packing in $T$ and therefore there exists a vertex (in fact, exactly one vertex), say $u$, in $N_T(v)\cap S'(T)$. Now, if $u'$ is the leaf adjacent to $u$, then $S=(S'(T) \setminus \{u\}) \cup \{u',v'\}$ is a $2$-packing in $T$, and so $\rho(T) \ge |S| > |S'(T)| = \rho(T)$, a contradiction.~\smallqed

\smallskip
\begin{claim} \label{claim2}
$d_T(x,S'(T) \setminus \{x\})=3$ for each $x\in S'(T)$.
\end{claim}
\noindent \emph{Proof of Claim~\ref{claim2}}. Since $S'(T)$ is a $2$-packing in $T$, $d_T(x,y)\ge 3$
for each pair of vertices $x, y\in S'(T)$, and therefore $d_T(x,S'(T) \setminus \{x\})\ge 3$ for each $x\in S'(T)$. Suppose that $d_T(v,S'(T) \setminus  \{v\})=k\ge 4$ for some vertex $v\in S'(T)$. Let $u\in S'(T)$ be a vertex such that $d_T(u,v)= d_T(v,S'(T)  \setminus \{v\})=k$, and let $v_0v_1 \ldots v_k$ be the $(v,u)$-path in $T$ where $v = v_0$ and $u = v_k$. If $k \ge 5$ and if $u'$ and $v'$ be the leaves adjacent to $u$ and $v$, respectively, then the set $S = (S'(T) \setminus \{u,v\}) \cup \{v', u', v_2\}$ is a $2$-packing in $T$, and so $\rho(T) \ge |S| > |S'(T)| = \rho(T)$, a contradiction. Hence, $k=4$. In this case, let $A = \{w\in S'(T)\colon d_T(v,w)=4\}$. Since $u \in A$, we note that $|A| \ge 1$. For each $w \in A\cup \{v\}$, let $w'$ be the (unique) leaf neighbor of $w$ in $T$. Now, the set
\[
S = (S'(T) \setminus (A \cup \{v\})) \cup (\{v',v_2\} \cup \{w'\colon w\in A\})
\]
is a $2$-packing in $T$, and so $\rho(T) \ge |S| > |S'(T)| = \rho(T)$, a contradiction.~\smallqed

\smallskip
\begin{claim}
\label{claim3}
If $x$ and $y$ are distinct vertices in $S'(T)$, then there is exactly one sequence of vertices, say $(x_0, x_1, \ldots, x_k)$ where $x = x_0$ and $y = x_k$, of distinct vertices in $S'(T)$ such that $d_T(x_{i-1},x_i)=3$ for $i \in [k]$.
\end{claim}

Before we present a proof of Claim~\ref{claim3}, as an illustration of the claim, if $T = G({\cal P})$ is the tree shown in Fig.~\ref{rys-do-podrodziny}(b), then the set $S'(T) = \{a,b,c,d,e,f\}$. Moreover, if $x = a$ and $y = f$, for example, then $x$ and $y$ are distinct vertices in $S'(T)$ and there is exactly one sequence $(x_0, x_1, x_2, x_3)$ where $x_0 = x$, $x_1=e$, $x_2=c$, and $x_3 = y$ of distinct vertices in $S'(T)$ such that $d_T(x_{i-1},x_i)=3$ for $i \in [3]$.

As a further illustration of Claim~\ref{claim3}, consider the tree $T$ shown in Fig.~\ref{rys-do-Claim 3}. If $S'(T)$ is the set of weak support vertices in $T$ and if $x = a$ and $y = x_6$, for example, then $x$ and $y$ are distinct vertices in $S'(T)$ and there is exactly one sequence $(x_0, x_1, x_2, x_3,x_4,x_5,x_6)$ where $x_0 = x$ and $x_6 = y$ of distinct vertices in $S'(T)$ such that $d_T(x_{i-1},x_i)=3$ for $i \in [6]$.

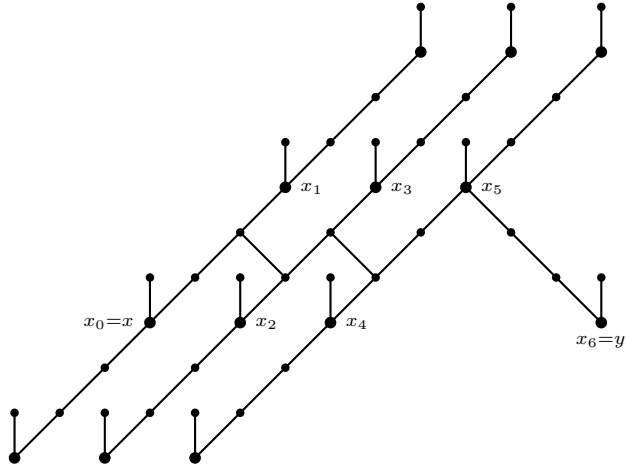
\begin{figure}[htb]
\begin{center} 
\begin{tikzpicture}[scale=0.6,style=thick,x=0.1cm,y=0.1cm] \def\vr{3pt}\def\vrr{2pt}
\draw (0,10)--(0,0)--(90,90)--(90,100);
\draw (20,10)--(20,0)--(110,90)--(110,100);
\draw (40,10)--(40,0)--(130,90)--(130,100);
\draw (30,30)--(30,40); \draw (60,60)--(60,70);
\draw (50,30)--(50,40); \draw (80,60)--(80,70);
\draw (70,30)--(70,40); \draw (100,60)--(100,70);
\draw (50,50)--(60,40); \draw (70,50)--(80,40);
\draw (100,60)--(130,30);
\draw (0,0) [fill=black] circle (\vr);
\draw (30,30) [fill=black] circle (\vr);
\draw (60,60) [fill=black] circle (\vr);
\draw (90,90) [fill=black] circle (\vr);
\draw (20,0) [fill=black] circle (\vr);
\draw (50,30) [fill=black] circle (\vr);
\draw (80,60) [fill=black] circle (\vr);
\draw (110,90) [fill=black] circle (\vr);
\draw (40,0) [fill=black] circle (\vr);
\draw (70,30) [fill=black] circle (\vr);
\draw (100,60) [fill=black] circle (\vr);
\draw (130,90) [fill=black] circle (\vr);
\draw (130,30) [fill=black] circle (\vr);
\draw (10,10) [fill=black] circle (\vrr);
\draw (20,20) [fill=black] circle (\vrr);
\draw (40,40) [fill=black] circle (\vrr);
\draw (50,50) [fill=black] circle (\vrr);
\draw (70,70) [fill=black] circle (\vrr);
\draw (80,80) [fill=black] circle (\vrr);
\draw (30,10) [fill=black] circle (\vrr);
\draw (40,20) [fill=black] circle (\vrr);
\draw (60,40) [fill=black] circle (\vrr);
\draw (70,50) [fill=black] circle (\vrr);
\draw (90,70) [fill=black] circle (\vrr);
\draw (100,80) [fill=black] circle (\vrr);
\draw (50,10) [fill=black] circle (\vrr);
\draw (60,20) [fill=black] circle (\vrr);
\draw (80,40) [fill=black] circle (\vrr);
\draw (90,50) [fill=black] circle (\vrr);
\draw (110,70) [fill=black] circle (\vrr);
\draw (120,80) [fill=black] circle (\vrr);
\draw (0,10) [fill=black] circle (\vrr);
\draw (30,40) [fill=black] circle (\vrr);
\draw (60,70) [fill=black] circle (\vrr);
\draw (90,100) [fill=black] circle (\vrr);
\draw (20,10) [fill=black] circle (\vrr);
\draw (50,40) [fill=black] circle (\vrr);
\draw (80,70) [fill=black] circle (\vrr);
\draw (110,100) [fill=black] circle (\vrr);
\draw (40,10) [fill=black] circle (\vrr);
\draw (70,40) [fill=black] circle (\vrr);
\draw (100,70) [fill=black] circle (\vrr);
\draw (130,100) [fill=black] circle (\vrr);
\draw (110,50) [fill=black] circle (\vrr);
\draw (120,40) [fill=black] circle (\vrr);
\draw (130,40) [fill=black] circle (\vrr);
\draw (130,30)--(130,40);
\draw (21,29.0) node {\small ${}^{x_0=x}$};
\draw (56,29.0) node {\small ${}^{x_2}$};
\draw (76,29.0) node {\small ${}^{x_4}$};
\draw (66,59.0) node {\small ${}^{x_1}$};
\draw (86,59.0) node {\small ${}^{x_3}$};
\draw (106,59.0) node {\small ${}^{x_5}$};
\draw (130,25.0) node {\small ${}^{x_6=y}$};
\end{tikzpicture}
\caption{A tree $T$ illustrating Claim~\ref{claim3}} \label{rys-do-Claim 3}
\end{center}
\end{figure}

\noindent \emph{Proof of Claim~\ref{claim3}}. It follows from Claim \ref{claim2} that $d_T(x,y)\ge 3$. To prove the desired result, we proceed by induction on $d_T(x,y)$. If $d_T(x,y)= 3$, then the sequence $(x,y)$ has the desired property and this establishes the base case. Assume, then, that the result holds for all pairs $x', y'\in S'(T)$ such that $3\le d_T(x',y') < q$, where $q \ge 4$. Assume that $x, y\in S'(T)$, $d_T(x,y)=q$, and let $y_0y_1 \ldots y_q$ be the $(x,y)$-path in $T$ where $x = y_0$ and $y = y_q$. Thus, $y_0 \in S'(T)$, and $y_1, y_2 \notin S'(T)$ (as $S'(T)$ is a $2$-packing).

We now prove that a vertex (and then exactly one vertex) belonging to $N_T(y_2) \setminus \{y_1\}$ is in $S'(T)$. Suppose, to the contrary, that $(N_T(y_2) \setminus \{y_1\}) \cap S'(T)= \emptyset$. Thus, $N_T[y_2]\cap S'(T)= \emptyset$. In this case, the set $S= (N_T(S'(T))\cap L(T))\cup \{y_2\}=L(T)\cup \{y_2\}$ is a $2$-packing in $T$, and so $\rho(T) \ge |S| = |S'(T)| + 1>|S'(T)| = \rho(T)$, a contradiction. Therefore, $(N_T(y_2) \setminus \{y_1\})\cap S'(T)\ne \emptyset$.

If $y_3 \in S'(T)$, then $3 \le d_T(y_3,y) < q$ and, by our inductive hypothesis, there is exactly one sequence $(z_0, z_1, \ldots, z_k)$ of distinct vertices in $S'(T)$ where $z_0 = y_3$ and $y = z_k$ and such that $d_T(z_{i-1},z_i)=3$ for $i \in [k]$. Thus, $(x_0, x_1, \ldots, x_{k+1})$ where $x = x_0$, $y = x_{k+1}$ and where $x_i = z_{i+1}$ for $i \in [k]$, is the desired sequence. Hence, we may assume that $y_3\notin S'(T)$, for otherwise the desired result follows.

Let $y_2'$ be the (unique) element of $(N_T(y_2) \setminus \{y_1,y_3\}) \cap S'(T)$. Thus, $3 \le d_T(y_2',y)<n$ and the inductive hypothesis guarantees that there is exactly one sequence $(z_0, z_1, \ldots, z_k)$ of distinct vertices in $S'(T)$ where $z_0 = y_2'$ and $z_k = y$ and such that $d_T(z_{i-1},z_i)=3$ for $i \in [k]$. Thus, $(x_0, x_1, \ldots, x_{k+1})$ where $x = x_0$ and where $x_{i+1} = z_{i}$ for $i \in [k]$, yielding the  desired sequence.~\smallqed

\smallskip
We now return to the proof of Theorem~\ref{obserwacja-O-prim-$2$-packing}, and are ready to prove that $T$ is in the family~$\cO'$. We shall prove that $T$ is isomorphic to $R({\cal P})$ for some tree $R$ and some family ${\cal P}$ of partitions of the sets $N_R(u)$ where $u\in V(R)$. Let $R=(V(R),E(R))$ be a graph with vertex set $V(R)= S'(T)$ and edge set $E(R)= \{uv\colon u, v\in V(R) \mbox{ and } d_T(u,v)=3\}$. It follows from Claim~\ref{claim3} that $R$ is a tree. Let ${\cal P}= \{{\cal P}(u)\colon u\in V(R)\}$ be a family, where for a vertex $u\in V(R)$, ${\cal P}(u)$ is the family $\{A_{ux}\colon x\in N_T(u) \setminus L(T)\}$ in which
\[
A_{ux} = \{y\in S'(T)\colon d_T(x,y)=2\}= S'(T)\cap N_T^3(u)\cap N_T^2(x)
\]
for $x\in N_T(u) \setminus L(T)$. We note that if $u\in V(R)$, then
\[
\begin{array}{lcl}
\displaystyle{  \bigcup_{x\in N_T(u) \setminus L(T)}A_{ux} }
& = & \displaystyle{ \bigcup_{x\in N_T(u) \setminus L(T)} (S'(T)\cap N_T^3(u)\cap N_T^2(x)) } \2 \\
& = & \displaystyle{ S'(T)\cap N_T^3(u) \cap \left( \bigcup_{x\in N_T(u) \setminus L(T)} N_T^2(x) \right) } \2 \\
& = & \displaystyle{ S'(T)\cap N_T^3(u) } \2 \\
& = & \displaystyle{ N_R(u) }.
\end{array}
\]

In addition, since
\[
A_{ux}\cap A_{uy} =S'(T)\cap N_T^3(u) \cap N_T^2(x)\cap N_T^2(y)
\]
and $T$ is a tree, the sets $A_{ux}$ and $A_{uy}$ are disjoint if $x, y\in N_T(u) \setminus L(T)$ and $x\ne y$. The above implies that ${\cal P}(u) =\{A_{ux}\colon x\in N_T(u) \setminus L(T)\}$ is a partition of the set $N_R(u)$ for $u\in V(R)$, and so
\[
\{u\}\times \cP(u) = \{(u,A_{ux})\colon x\in N_T(u) \setminus L(T)\}
\]
for every vertex $u\in V(R)$. Let us consider the graph $R(\cP)$ for the above defined set $R$ and the partition $\cP$. By construction, $R(\cP)$ is a graph with vertex set
\[
V(R(\cP)) = V(R)\cup (V(R)\times\{1\}) \cup \left( \bigcup_{u\in V(R)}\{(u,A_{ux})\colon x\in N_T(u) \setminus L(T)\} \right)
\]
and edge set $E(R(\cP)) = E_1 \cup E_2 \cup E_3$ where
\[
\begin{array}{lcl}
E_1 & = & \displaystyle{ \{u(u,1)\colon u\in V(R)\}  }, \2 \\
E_2 & = & \displaystyle{ \bigcup_{u\in V(R)}\{u(u,A)\colon A\in \cP(u)\} }, \2 \\
E_3 & = & \displaystyle{ \bigcup_{uv\in E(R)}\{(u,A)(v,B)\colon  A\in \cP(u), B\in \cP(v), u\in B, v\in A\} }. \\
\end{array}
\]
Equivalently, the edge sets $E_2$ and $E_3$ are the sets
\[
E_2 = \bigcup_{u\in V(R)}\{u(u,A_{ux})\colon x\in N_T(u) \setminus L(T)\}
\]
and
\[
E_3 = \bigcup_{uv\in E(R)}\{(u,A_{ux})(v,B_{vy})\colon x\in N_T(u) \setminus L(T), y\in N_T(v) \setminus L(T), u\in B_{vy}, v\in A_{ux}\}.
\]

We now let
\[
\varphi \colon V(R(\cP)) \to V(T)
\]
be a function defined in such a way that $\varphi(x)=x$ for all $x\in V(R)$, $\varphi((x,1))= l_x$ if $(x,1)\in V(R)\times \{1\}$ and $l_x$ is the only leaf adjacent to $x$ in $T$, and, finally, $\varphi((u,A_{ux}))= x$ if $u\in V(R)$, $x\in N_T(u) \setminus L(T)$, and $A_{ux}\in \cP(u)$. It is immediate from the definitions of $R$, $\cP$, and $R(\cP)$ that the function $\varphi$ is an isomorphism of $T$ and $R(\cP)$. This completes the proof.
\end{proof}

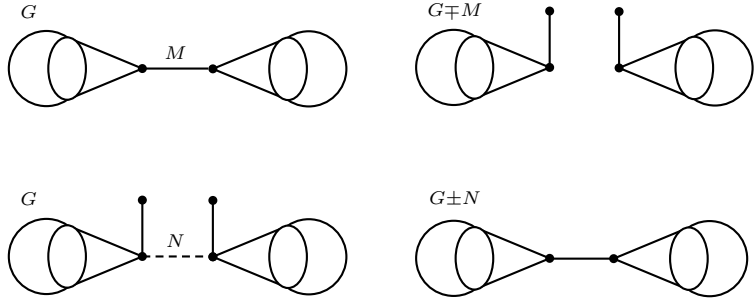
\begin{figure}[htb] \begin{center} 
\begin{tikzpicture}[scale=0.5,style=thick,x=0.1cm,y=0.1cm] \def\vr{2.6pt}
\draw (-40,50) node {\begin{tikzpicture}[scale=0.5,style=thick,x=0.1cm,y=0.1cm]
\draw (-20,28.3)--(0,20)--(-20,11.7);
\draw (0,20)--(17.5,20);
\draw (0,20) [fill=black] circle (\vr);
\draw[fill=white] (-25.5,20) ellipse (10 and 10);
\draw[fill=white] (-20,20) ellipse (5 and 8.3);
\draw (-30,34) node {\small ${}^G$};
\draw (9,23) node {\small ${}^M$};
\end{tikzpicture}};
\draw (5,45.3) node {\begin{tikzpicture}[scale=0.5,style=thick,x=0.1cm,y=0.1cm] \def\vr{2.6pt}
\draw (20,28.3)--(0,20)--(20,11.7);
\draw (0,20) [fill=black] circle (\vr);
\draw[fill=white] (25.5,20) ellipse (10 and 10);
\draw[fill=white] (20,20) ellipse (5 and 8.3);
\end{tikzpicture}};
\draw (60,50.3) node {\begin{tikzpicture}[scale=0.5,style=thick,x=0.1cm,y=0.1cm] \def\vr{2.6pt}
\draw (-20,28.3)--(0,20)--(-20,11.7);
\draw (0,20)--(0,35);
\draw (0,20) [fill=black] circle (\vr);
\draw (0,35) [fill=black] circle (\vr);
\draw[fill=white] (-25.5,20) ellipse (10 and 10);
\draw[fill=white] (-20,20) ellipse (5 and 8.3);
\draw (-25,34) node {\small ${}^{G\mp M}$};
\end{tikzpicture}};
\draw (113,48.5) node {\begin{tikzpicture}[scale=0.5,style=thick,x=0.1cm,y=0.1cm] \def\vr{2.6pt}
\draw (20,28.3)--(0,20)--(20,11.7);
\draw (0,20) [fill=black] circle (\vr);
\draw (0,35) [fill=black] circle (\vr);
\draw (0,20)--(0,35);
\draw[fill=white] (25.5,20) ellipse (10 and 10);
\draw[fill=white] (20,20) ellipse (5 and 8.3);
\end{tikzpicture}};
\draw (-40,0) node {\begin{tikzpicture}[scale=0.5,style=thick,x=0.1cm,y=0.1cm] \def\vr{2.6pt}
\draw (-20,28.3)--(0,20)--(-20,11.7);
\draw (0,20) [fill=black] circle (\vr);
\draw (0,35) [fill=black] circle (\vr);
\draw (0,20)--(0,35);
\draw[densely dashed] (0,20)--(17.5,20);
\draw[fill=white] (-25.5,20) ellipse (10 and 10);
\draw[fill=white] (-20,20) ellipse (5 and 8.3);
\draw (-30,34) node {\small ${}^G$};
\draw (9,23) node {\small ${}^N$};
\end{tikzpicture}};
\draw (5,-1.75) node {\begin{tikzpicture}[scale=0.5,style=thick,x=0.1cm,y=0.1cm] \def\vr{2.6pt}
\draw (20,28.3)--(0,20)--(20,11.7);
\draw (0,20) [fill=black] circle (\vr);
\draw (0,20) [fill=black] circle (\vr);
\draw (0,35) [fill=black] circle (\vr);
\draw (0,20)--(0,35);
\draw[fill=white] (25.5,20) ellipse (10 and 10);
\draw[fill=white] (20,20) ellipse (5 and 8.3);
\end{tikzpicture}};
\draw (60,0) node {\begin{tikzpicture}[scale=0.5,style=thick,x=0.1cm,y=0.1cm] \def\vr{2.6pt}
\draw (-20,28.3)--(0,20)--(-20,11.7);
\draw (0,20) [fill=black] circle (\vr);
\draw[fill=white] (-25.5,20) ellipse (10 and 10);
\draw[fill=white] (-20,20) ellipse (5 and 8.3);
\draw (-25,35) node {\small ${}^{G\pm N}$};
\end{tikzpicture}};
\draw (103,-5.2) node {\begin{tikzpicture}[scale=0.5,style=thick,x=0.1cm,y=0.1cm] \def\vr{2.6pt}
\draw (20,28.3)--(0,20)--(20,11.7);
\draw (0,20)--(-18,20);
\draw (0,20) [fill=black] circle (\vr);
\draw[fill=white] (25.5,20) ellipse (10 and 10);
\draw[fill=white] (20,20) ellipse (5 and 8.3);
\end{tikzpicture}};
\end{tikzpicture} \end{center} \vskip -0.4 cm
\caption{$M$ and $N$ are matchings in $G$ and  $\overline{G}$, respectively}\label{tree-operacja22222}
\end{figure}

\smallskip
In this final subsection, we give a simple construction that makes it possible to build  trees belonging to the family $\cO$ from smaller trees belonging to the family $\cO'$. Let $M$ be a matching in a graph $G$. By $G\mp M$ we denote the graph obtained from $G-M$ by adding exactly one new pendant edge to each vertex covered by $M$. On the other hand, let $N$ be a matching in the complement $\overline{G}$ of $G$, and assume that every vertex covered by $N$ is a weak support in $G$. Then by $G\pm N$ we denote the graph obtained from $G\cup N$ by removing the leaf neighbor of every vertex in $N$. It follows from these definitions that if $M$ is a matching in $G$, then $(G\mp M)\pm M=G$. Similarly, if $N$ is a matching in $\overline{G}$ and every vertex covered by $N$ is a weak support in $G$, then $(G\pm N)\mp N=G$, see Fig.~\ref{tree-operacja22222} for an illustration.

As a consequence of Theorem~\ref{obserwacja-O-prim-$2$-packing}, we have the following corollary. 

\begin{cor} 
\label{cor:1}
Let $u$ and $v$ be weak support vertices in disjoint trees $T_1$ and $T_2$, respectively, and let $T_{uv}$ be the tree obtained from the union $T_1 \cup T_2$ by adding the edge $uv$ and removing the leaf neighbor adjacent to $u$ and $v$, respectively. Then the tree $T_{uv}$ is in Class~$3$ if both $T_1$ and $T_2$ are in Class~$3$.
\end{cor}
\begin{proof} Assume that $T_1$ and $T_2$ are in Class~$3$, and let $F_1$ and $F_2$ be the subsets of $E(T_1)$ and $E(T_2)$ having properties (4a)--(4c) of Theorem \ref{new-characterization} in $T_1$ and $T_2$, respectively. Let $uu'$ and $vv'$ be pendant edges incident with $u$ and $v$ in $T_1$ and $T_2$, respectively. Then the set $(F_1 \cup F_2 \cup \{uv\}) \setminus \{uu', vv'\}$ has properties (4a)--(4c) of Theorem~\ref{new-characterization} in $T_{uv}$ and, therefore, $T_{uv}$ is in Class~$3$.
\end{proof}

We remark that the example of a path $P_{10}$ illustrates that the converse implication of Corollary~\ref{cor:1} does not hold without additional assumptions on the edge $uv$ and vertices $u$ and $v$ in $T_{uv}$.

\begin{thm} 
If $T$ is a tree of order at least~$6$, then  $T\in \cO$ if and only if
$T\in \cO'$ or every component of $T\mp M$ belongs to $\cO'$ for some matching $M$ in $T$. 
\end{thm}
\begin{proof} 
Assume that $T\in \cO$. Then in $T$ there is a uniquely determined subset $F$ of $E(T)$ that has properties (4a)--(4c) of Theorem~\ref{new-characterization}. If every edge in $F$ is a~pendant edge in $T$, then it follows from (4b) that the set $S'(T)$ (of weak supports of $T$) is a $2$-packing in $T$.

We claim that $S'(T)$ is a maximum $2$-packing in $T$. Suppose, to the contrary, that there is a~$2$-packing $S$ in $T$ such that $|S|>|S'(T)|$. Since $|N_T[x] \cap S|\le 1$ for each $x\in S'(T)$, the inequality $|S|>|S'(T)|$ implies that the set $S \setminus N_T[S'(T)]$ is nonempty, say $x_0\in S \setminus N_T[S'(T)]$. Thus, $d_T(x_0,S'(T))\ge 2$ and, therefore, if $e$ and $f$ are pendant edges in $T$ belonging to distinct components of $T-x_0$, then no $e$\,--$f$ sequence has property (4c), a~contradiction. This proves that $S'(T)$ is a maximum $2$-packing in $T$ and $T\in \cO'$, by Theorem~\ref{obserwacja-O-prim-$2$-packing}. 

We may therefore assume that not every edge belonging to $F$ is a~pendant edge in $T$. Let $M$ be the set of non-pendant edges in $T$ belonging to $F$. Now the graph $T\mp M$ is disconnected, it has $k=|M|+1$ components, say $T_1, \ldots, T_k$, and $S'(T)\cup V(M)$ is the set of weak support vertices in $T \mp M$ (where $V(M)$ is the set of vertices covered by $M$). From the fact that $F$ has properties (4a)--(4c) of Theorem \ref{new-characterization} it follows that the set of pendant edges in $T_i$ where $i \in [k]$ has properties (4a)--(4c) of Theorem~\ref{new-characterization}. Then, as in the beginning of this proof, the set $S'(T_i)$ of weak support vertices of $T_i$ is a maximum $2$-packing in $T_i$, and, consequently, $T_i\in \cO'$, by Theorem~\ref{obserwacja-O-prim-$2$-packing} for $i \in [k]$.

Assume now that $M$ is a matching in $T$ such that every component, say $T_1, \ldots, T_k$, of $T\mp M$ belongs to $\cO'$. Then the set of weak support vertices of $T_i$ is a maximum $2$-packing in $T_i$ for $i \in [k]$ (and the set of weak support vertices of $T\mp M$ is a~maximum $2$-packing in  $T\mp M$). Consequently, the set $F_i$ of pendant edges of $T_i$ has properties (4a)--(4c) of Theorem~\ref{new-characterization} in $T_i$. This implies that the set of edges $M \cup \bigcup_{i=1}^k (F_i \cap E(T))$  has properties (4a)--(4c) of Theorem~\ref{new-characterization} in $T$. From these observations and from Theorem~\ref{new-characterization} we infer that $T\in \cO$.
\end{proof}


\medskip
\begin{thebibliography}{99}


\bibitem{Favaron-Karami-Khoeilar-Sheikholeslami-1} O. Favaron, H. Karami, R. Khoeilar,  and S. M. Sheikholeslami, A new bound on the total domination subdivision number. \textit{Graphs Combin.} \textbf{25} (2009), 41--47.

\bibitem{Favaron-Karami-Khoeilar-Sheikholeslami-2} O. Favaron, H. Karami, R. Khoeilar,  and S. M. Sheikholeslami, On the total domination subdivision number in some classes of graphs. \textit{J. Combin. Optim.}  \textbf{20} (2010), 76--84.

\bibitem{Favaron-Karami-Khoeilar-Sheikholeslami-3} O. Favaron, H. Karami, R. Khoeilar, and S. M. Sheikholeslami, Matching and total domination subdivision number of graphs with few $C_4$. \textit{Discuss. Math. Graph Theory} \textbf{30} (2010), 611--618.

\bibitem{Favaron-Karami-Khoeilar-Sheikholeslami-4} O. Favaron, H. Karami, R. Khoeilar,  and S. M. Sheikholeslami, On the total domination subdivision number in graphs. \textit{Bull. Malays. Math. Sci. Soc.} \textbf{37} (2014), 173--180.

\bibitem{Favaron-Karami-Sheikholeslami-1} O. Favaron, H. Karami,  and S. M. Sheikholeslami, Total domination and total domination subdivision numbers. \textit{Australas. J. Combin.} \textbf{38} (2007), 229--235.

\bibitem{Favaron-Karami-Sheikholeslami-2} O. Favaron, H. Karami,  and S. M. Sheikholeslami, Total domination and total domination subdivision number of a graph and its complement. \textit{Discrete Math.} \textbf{308} (2008), 4018--4023.

\bibitem{HaHeHe-20} T. W. Haynes, S. T. Hedetniemi, and M. A. Henning  (eds), \emph{Topics in Domination in Graphs}. Series: Developments in Mathematics, Vol. 64, Springer, Cham, 2020. viii + 545 pp.

\bibitem{HaHeHe-21} T. W. Haynes, S. T. Hedetniemi, and M. A. Henning  (eds), \emph{Structures of Domination in Graphs}. Series: Developments in Mathematics, Vol. 66, Springer, Cham, 2021. viii + 536 pp.

\bibitem{HaHeHe-23} T. W. Haynes, S. T. Hedetniemi, and M. A. Henning, \emph{Domination in Graphs: Core Concepts}. Series: Springer Monographs in Mathematics, Springer, Cham, 2023. xx + 644 pp.

\bibitem{Haynes-Hedetniemi-van der Merwe} T. W. Haynes, S. T. Hedetniemi, and L. C. van der Merwe, Total domination subdivision numbers. \textit{J. Combin. Math. Combin. Comput.} \textbf{44} (2003), 115--128.

\bibitem{Haynes-Henning-HopkinsDMGT} T. W. Haynes, M. A. Henning,  and L. Hopkins, Total domination subdivision numbers of graphs. \textit{Discuss. Math. Graph Theory} \textbf{24} (2004), 457--467.

\bibitem{Haynes-Henning-Hopkins} T. W. Haynes, M. A. Henning,  and  L. Hopkins, Total domination subdivision numbers of trees. \textit{Discrete Math.} \textbf{286} (2004), 195--202.

\bibitem{HeYe-13} M. A. Henning and A. Yeo, \emph{Total Domination in Graphs.} Series: Springer Monographs in Mathematics, Springer, Cham, 2013. xiv + 178 pp.

\bibitem{Jeyanthi-Hemalatha-Davvaz} P. Jeyanthi, G. Hemalatha,  and B. Davvaz, Total domination subdivision number in strong product graph. \textit{Amer. J. Appl. Math. Stat.} \textbf{2} (2014), 216--219.

\bibitem{Karami-Khodkar-Khoeilar-Sheikholeslami} H. Karami, A. Khodkar, R. Khoeilar,  and S. M. Sheikholeslami, Trees whose total domination subdivision number is one. \textit{Bull. Inst. Combin. Appl.} \textbf{53} (2008), 57--67.

\bibitem{Karami-Khoeilar-Sheikholeslami} H. Karami, R. Khoeilar,  and S. M. Sheikholeslami,
    The total domination subdivision number in graphs with no induced $3$-cycle and $5$-cycle. \textit{J. Combin. Optim.} \textbf{25} (2013), 91--98.

\bibitem{Karami-Khodkar-Sheikholeslami} H. Karami, A. Khodkar,  and S. M. Sheikholeslami, An upper bound for total domination subdivision numbers of graphs. \textit{Ars Combin.}  \textbf{102} (2011), 321--331.

\bibitem{Khoeilar-Karami-Sheikholeslami} R. Khoeilar, H. Karami,  and S. M. Sheikholeslami,
    On two conjectures concerning total domination subdivision number in graphs. \textit{J. Combin. Optim.} 38 (2019), 333--340.

\bibitem{Sheikholeslami} S. M. Sheikholeslami, On the total domination subdivision numbers in graphs. \textit{Cent. Eur. J. Math.} \textbf{8} (2010), 468--473.

\bibitem{Soltankhah} N. Soltankhah, On the total domination subdivision numbers of grid graphs. \textit{Int. J. Contemp. Math. Sciences} \textbf{5} (2010), 2419--2432.

\end{thebibliography}
\end{document}